\documentstyle[12pt]{article}
\makeatletter\@addtoreset{equation}{section}\makeatother

\begin{document}
\begin{center}
{\large\bf THE DETERMINATION OF THE INTERNAL STRUCTURE OF THE SUN
BY THE DENSITY DISTRIBUTION}\\[0,5cm]
{\bf H.J. HAUBOLD}\\
UN Outer Space Office, Vienna International Centre, Vienna,
Austria\\[0,5cm]
and\\[0,5cm]
{\bf A.M. MATHAI}\\
Department of Mathematics and Statistics, McGill University,
Montreal, Canada\\[1,0cm]
{\bf ABSTRACT}
\end{center}
\bigskip
This paper examines a number of analytic models which can describe
the gravitationally stabilized fusion reactor of the Sun. An
analytical multiparameter model is shown to reproduce various
structure parameters such as density, mass, pressure and
temperature throughout the solar core, which have been obtained by
solving numerically the system of solar structure differential
equations. For simplicity, numerical and analytical results are
discussed in detail for a two-parameter model.
\section{INTRODUCTION}
The numerical approach to the problem of studying solar structure
is to go for numerical solutions of the underlying system of
differential equations (A. Noels, R. Papy, and F. Remy 1993). Even
for a simple main-sequence star in hydrostatic equilibrium, like
the Sun, at least four nonlinear differential equations are to be
dealt with to obtain a detailed picture of the run of physical
variables throughout the star. Another approach is to search for
analytical solutions of these differential equations under
justified physical simplifications (Haubold and Mathai 1992).\par
In this paper we start with numerical data given in Sears (1964)
for various analytical models which will fit the data on density in
the gravitationally stabilized solar core. Let r be the distance
from the center of the Sun to an arbitrary point in its interior,
R the solar radius, $\rho(r)$ the density at $r$, $\rho_c$ the
density at the center and $y=r/R$. For the range $0\leq y\leq 0.3$
consider the following models for the density $\rho(r)$.
\begin{eqnarray}
u & = & \frac{\rho(r)}{\rho_c} =
1-4.94y+6.67y^2-2.73y^3\label{eq:1.1}\\
u & = & 1-4y+2y^2+2y^3-y^4\label{eq:1.2}\\ 
u & = & (1-\sqrt{y})(1-y^3)^{64}\label{eq:1.3}\\ 
u & = & (1-y^{3/2})^{16}\label{eq:1.4}\\
u & = & (1-\sqrt{y})(1-y^3)^{64}(1-y)\label{eq:1.5}\\ 
u & = & (1-y^{1.48})^{14}\label{eq:1.6}\\
u & = & (1-y^{1.48})^{13}\label{eq:1.7}\\
u & = & (1-y^{1.28})^{10}\label{eq:1.8}
\end{eqnarray}
The following table gives $y$, the numerically obtained data on
$u=\rho(r)/\rho_c$ from Sear's (1964) and the estimated values
under models (1.1) to (1.8).\par
\bigskip
\noindent
Table 1 density distribution\par
\bigskip
\begin{tabbing}
\=Model(11) \=Mod(Sears,1964) \=Model0000 \=Model0000 \=Model0000
\=Model0000\= \kill
\>$y=\frac{r}{R}$ \>$\frac{\rho(r)}{\rho_c}$(Sears, 1964) \>Model
\>Model \>Model \>Model \\
\> \> \>(1.1) \>(1.2) \>(1.3) \>(1.4) \\ [0.5cm]
\>0.0864 \>0.6519 \>0.6213 \>0.6720 \>0.6795 \>0.6626\\
\>0.1153 \>0.5253 \>0.5149 \>0.5690 \>0.5987 \>0.5283\\
\>0.1441 \>0.3856 \>0.4185 \>0.4710 \>0.5722 \>0.4065\\
\>0.1873 \>0.2810 \>0.2908 \>0.3340 \>0.3720 \>0.2588\\
\>0.2161 \>0.1994 \>0.2164 \>0.2490 \>0.2796 \>0.1837\\
\>0.2450 \>0.1424 \>0.1499 \>0.1700 \>0.1957 \>0.1264\\
\>0.2882 \>0.0962 \>0.0649 \>0.0560 \>0.0985 \>0.0679 
\end{tabbing}
\clearpage
\noindent
Table 1 continued
\begin{tabbing}
\=Model (1.5) \=Model (1.6) \=Model (1.7) \=Model (1.8) \= \kill\\ 
\>Model \> Model \>Model \> Model \\
\>(1.5) \>(1.6) \>(1.7) \>(1.8)\\ [0,5cm]
\>0.6208 \>0.6849 \>0.7037 \>0.6418\\
\>0.5297 \>0.5573 \>0.5811 \>0.5229\\
\>0.4435 \>0.4404 \>0.4669 \>0.4179\\
\>0.3023 \>0.2928 \>0.3196 \>0.2885\\
\>0.2192 \>0.2159 \>0.2409 \>0.2202\\
\>0.1478 \>0.1551 \>0.1772 \>0.1649\\
\>0.0701 \>0.0885 \>0.1053 \>0.1033
\end{tabbing}
Model (1.1) results from a least square fit of a third degree
polynomial to the data of Sears (1964). Model (1.2) comes from
successive eliminations. Even though these models can adequately
describe the data in the range $0\leq y \leq 0.3$ the equation u=0
has real roots in (0,1) besides the root at $y=1$. Hence these
cannot be used
as models for the density in [0,1]. From models (1.3) to (1.8) it
is
evident that a model of the type
\begin{equation}
u=(1-y^{a_1})^{b_1}(1-y^{a_2})^{b_2}\ldots(1-y^{a_k})^{b_k},
\label{eq:1.9}
\end{equation}
where $a_i>0, b_i=0,1,2,\ldots, i=1,\ldots,k$ can be an excellent
fit to the data. For simplicity we will study a model of the type
\begin{equation}
u=(1-y^\delta)^\gamma,\label{eq:1.10}
\end{equation}
where $\delta>0$ and $\gamma$ a positive integer. Illustration will
be given for the case $\delta=1.28$ and $\gamma=10.$
\section{A Solar Model}
We assume that rotation and magnetic fields have no impact on the
internal structure of the Sun and the gravitational force directed
inward and the gas pressure force directed outward keep the gas
sphere in hydrostatic equilibrium. The ratio of gas pressure and
radiation pressure increases towards the center of the Sun but
never exceeds $10^{-3}$. Thus radiation pressure can be neglected
for the purpose of considering hydrostatic equilibrium. Assume that
the density $\rho(r)$ varies from $\rho_c$ at $r=0$ to
zero at $r=R$ in the following fashion:
\begin{equation}
\rho(r)=\rho_c[1-(\frac{r}{R})^\delta]^\gamma,\; \delta>0,\;
\gamma\;
\mbox{a positive integer.}
\end{equation}
Then the distribution of mass in terms of radius is given by
$$\frac{dM(r)}{dr}=4\pi r^2\rho (r),$$
so that
\begin{eqnarray}
M(r) & = & 4\pi \int^r_0dtt^2\rho(t)\nonumber \\ 
& = & 4\pi \rho_c \int^r_0dtt^2[1-(\frac{t}{R})^\delta]^\gamma
\nonumber \\
& = & 4\pi \rho_c\sum^\gamma_{m=0}(-1)^m(^\gamma_m)R^3
\int^{r/R}_0 dz z^{m\delta +2} \nonumber \\
& = &\frac{4\pi \rho_c}{3}R^3(\frac{r}{R})^3\;_2F_1(-\gamma,
\frac{3}{\delta};\frac{3}{\delta}+1;(\frac{r}{R})^\delta),
\end{eqnarray}
where $_2F_1$ is a Gauss' hypergeometric function, see for example,
Mathai (1993). M(r) gives the total mass at $r=R$, that is,
\begin{eqnarray}
M(R) & = & \frac{4\pi \rho_c R^3}{3}\; _2F_1 (-\gamma,
\frac{3}{\delta}; \frac{3}{\delta}+1; 1) \nonumber \\
& = & \frac{4\pi \rho_c
R^3}{3}\frac{\gamma!}{(\frac{3}{\delta}+1)(\frac{3}{\delta}+2)
\ldots
(\frac{3}{\delta}+\gamma)},
\end{eqnarray}
which is evaluated by using the formula
\begin{equation}
_2F_1(a,b;c;1)=\frac{\Gamma(c)\Gamma(c-a-b)}{\Gamma(c-a)\Gamma(c-
b)}.
\end{equation}
From (2.2) and (2.3) we can get a formula for estimating $\rho_c$
also. That is,
\begin{equation}
\rho_c =
\frac{3M(R)}{4\pi R^3}\frac{(\frac{3}{\delta}+1)(\frac{3}{\delta}+
2)\ldots(\frac{3}{\delta}+\gamma)}{\gamma!}.
\end{equation}
The numerical values for the Sun from Sears (1964) are
\begin{eqnarray}
\rho_c & = & 158 \mbox{gcm}^{-3}\nonumber \\
M(R) & = & 1.991 \times 10^{33} g \label{line2}\\
R & = & 6.96 \times 10^{10}cm. \nonumber
\end{eqnarray}
Substituting these in (2.5) we have
\begin{equation}
\frac{(\frac{3}{\delta}+1)(\frac{3}{\delta}+2)\ldots(\frac{3}
{\delta}+\gamma)}{\gamma!} = 112.08.
\end{equation}
Hence a pair of $(\delta, \gamma)$ which satisfy (2.7) is a good
choice. We have many choices for $(\delta, \gamma)$ even if
$\gamma$ is kept as a positive integer. Keeping $\gamma$ a positive
integer, for simplicity, it is found that $(\delta=1.28,
\gamma=10)$ is a
convenient choice in the sense that these values give good
estimates for $M(R)$ and $R$ also, close to the realistic values
and the estimates for the density are not far off.\par
Note that (2.7) is obtained by using the readings in (2.6). From
(2.7) one can get a set of positive values for $\delta$ and
$\gamma$. These are given here for some integer values of
$\gamma$:\\
Table 2 $\gamma$ and $\delta$
\begin{tabbing}
\=00000000 \=00000000 \=00000000 \=00000000 \=00000000 \=00000000
\=00000000 \= \kill \\
\>$\gamma$ \>2 \>3 \>4 \>5 \>6 \>7 \\
\>$\delta$ \>0.2225 \>0.4412 \>0.6265 \>0.7802 \>0.9098 \>1.0211 \\
[0.5cm]
\>\>8 \>9 \>10 \>11 \>12 \>13 \>14 \\
\>\>1.1182 \>1.2043 \>1.2814 \>1.3512 \>1.4149 \>1.4735 \>1.5276 \\
[0,5cm]
\>\>15 \>16 \>17 \>18 \>19 \>20 \\
\>\>1.5780 \>1.6251 \>1.6692 \>1.7108 \>1.75 \>1.7872
\end{tabbing}
From (2.2) and (2.3) the relative mass is given by
\begin{equation}
\frac{M(r)}{M(R)}=\frac{(\frac{3}{\delta}+1)(\frac{3}{\delta}+2)
\ldots(\frac{3}{\delta}+\gamma)}{\gamma!}(\frac{r}{R})^3\;_2F_1(-
\gamma, \frac{3}{\delta};\frac{3}{\delta}+1;(\frac{r}{R})^\delta).
\end{equation}
This is computed for $\delta=1.28$, $\gamma=10$ and compared with
Sear's (1964) numerical results in the following table 3.
\clearpage
\noindent
Table 3 mass distribution
\begin{tabbing}
\=0,000000000 \=(analytic modell)000 \=0,000000 \= \kill
\>$y=\frac{r}{R}$ \>$\frac{M(r)}{M(R)}$ \>$\frac{M(r)}{M(R)}$\\ 
\> \>(analytic model) \>(Sears 1964) \\ [0,5cm]
\>0.0864 \>0.0533 \>0.05\\
\>0.1153 \>0.1106 \>0.1\\
\>0.1441 \>0.1866 \>0.2\\
\>0.1873 \>0.3249 \>0.3\\
\>0.2161 \>0.4241 \>0.4 \\
\>0.2450 \>0.5225 \>0.5\\
\>0.2882 \>0.6576 \>0.6\\
\end{tabbing}
\par
Assuming that the pressure $P(r)$ at the center of the Sun is $P_c$
and that at the surface is zero we have
\begin{equation}
P(r)=P_c-G\int^r_0 dt\frac{M(t)\rho(t)}{t^2},
\end{equation}
where G is the gravitational constant. That is,
\begin{eqnarray}
P(r)& = & P_c-\frac{4\pi G}{\delta}\rho^2_c
R\sum^\gamma_{m=0}\frac{(-\gamma)_m}{m!}\int^r_0 dt
\frac{(\frac{t}{R})^{m\delta+1}}{(\frac{3}{\delta}+m)}[1-
(\frac{t}{R})^\delta]^\gamma \nonumber \\
& = & P_c-\frac{4\pi G}{\delta^2}
\rho^2_c R^2\sum^\gamma_{m=0}\frac{(-\gamma)_m(\frac{r}{R})^
{m\delta +2}}{(\frac{3}{\delta}+m)(\frac{2}{\delta}+m)}
\times \nonumber \\
& &_2F_1(-\gamma,\frac{2}{\delta}+m;\frac{2}{\delta}+m+1;
(\frac{r}{R})^\delta).\label{line3}
\end{eqnarray}
Assuming that the pressure at the surface is zero, that is
$P(R)=0$, we obtain
\begin{eqnarray}
P_c & = & \frac{4\pi G}{\delta^2}\rho^2_c R^2
\sum^\gamma_{m=0}\frac{(-
\gamma)_m}{m!}  
\frac{1}{(\frac{3}{\delta}+m)(\frac{2}{\delta}+m)}\times \nonumber
\\
& & \frac{\gamma!}{(\frac{2}{\delta}+m+1)(\frac{2}{\delta}+m+2)
\ldots(\frac{2}{\delta}+m+\gamma)}. \label{line2}
\end{eqnarray}
Then
\begin{eqnarray}
P(r)& = & \frac{4\pi
G}{\delta^2}\rho^2_c R^2\sum^\gamma_{m=0}\frac{(-
\gamma)_m}{m!}\frac{1}{(\frac{3}{\delta}+m)(\frac{2}{\delta}+m)}
\times \nonumber\\
& &\left[\frac{\gamma!}{(\frac{2}{\delta}+m+1)\ldots(\frac{2}
{\delta}+m+
\gamma)}-(\frac{r}{R})^{m\delta+2}\right.\times \nonumber \\
& &\left._2F_1(-\gamma,\frac{2}{\delta}+m;\frac{2}{\delta}+m+1;
(\frac{r}{R})^\delta)\right]
\end{eqnarray}
The above expansion can also be written as a double sum. In this
case
\begin{eqnarray}
P(r) & = & P_c-\frac{4\pi
G}{6}\rho^2_cr^2\sum^\gamma_{m=0}\sum^\gamma_{n=0}\frac{(-
\gamma)_m}{m!}\frac{(-\gamma)_n}{n!}\times \nonumber \\
& &
\frac{(\frac{2}{\delta})_m(\frac{3}{\delta})_m(\frac{2}{\delta})_
{m+n}}{(\frac{2}{\delta}+1)_m(\frac{3}{\delta}+1)_m(\frac{2}
{\delta}+1)_{m+n}}\left[(\frac{r}{R})^\delta\right]^m\left[(\frac
{r}{R})^
\delta\right]^n.\label{line2}
\end{eqnarray}
This double sum can be written in terms of a Kamp\'{e} de
F\'{e}riet's
function. That is,
\begin{equation}
P(r)=P_c-\frac{2}{3} \pi G \rho_c^2r^2
F^{1:3:1}_{1:2:0}\left[\left(^{(\frac{r}{R})^\delta}_{(\frac{
r}{R})^\delta}\right)\left|^
{\frac{2}{\delta}:-\gamma, \frac{3}{\delta}, \frac{2}{\delta}:-
\gamma}_{\frac{2}{\delta}+1:\frac{3}{\delta}+1,\frac{2}{\delta}+1:
}\right]\right.
\end{equation}
and then assuming that P(R)=0,
\begin{equation}
P_c=\frac{2}{3}\pi G \rho_c^2R^2
F^{1:3:1}_{1:2:0}\left[(^1_1)\left|^{\frac{2}{\delta}:-\gamma,
\frac{3}{\delta}, \frac{2}{\delta}:-\gamma}_{\frac{2}{\delta}+1:
\frac{3}{\delta}+1, \frac{2}{\delta}+1:}\right]\right.,
\end{equation}
where the Kamp\'{e} de F\'{e}riet's series is defined by the
following:
\begin{eqnarray}
& F^{p:q:k}_{r:m:n}& \left[\left(^x_y\right)\mid^{(a_p):(b_q):
(c_k)}_{(\alpha_r):(\beta_m):
(\gamma_n)}\right] = \label{line1} \\
& & =\sum^\infty_{m'=0}\sum^\infty_{n'=0}
\frac{\left[\Pi^p_{j=1}(a_j)
_{m'+n'}
\right]\left[\Pi^q_{j=1}(b_j)_{m'}\right]\left[\Pi^k_{j=1}(c_j)_
{n'}\right]}{\left[\Pi^r_{j=1}(\alpha_j)
_{m'+n'}\right]\left[\Pi^m_{j=1}(\beta_j)_{m'}
\right]\left[\Pi^n
_{j=1}(\gamma_j)_{n'}\right]}
\frac{x^{m'}y^{n'}}{m'!n'!}.\nonumber
\end{eqnarray}
For a discussion of Kamp\'{e} de F\'{e}riet's function see
Srivastava and Karlsson (1985). Note that (2.14) is a  polynomial
since $\gamma$ is a positive integer and hence convergence
conditions do not arise in (2.14). The form in (2.10) is the most
appropriate for computational purposes and the computations given
later are done by using (2.10).\par
The temperature is given by the equation of state of the perfect
gas
\begin{equation}
T(r)=\frac{\mu}{kN_A}{P(r)}{\rho(r)},
\end{equation}
where $\mu$ is the mean molecular weight, k is Boltzmann's constant
and $N_A$ is Avogadro's number.
Let
\begin{eqnarray}
g(r)& = & \frac{1}{\delta^2}\sum^\gamma_{m=0}\frac{(-
\gamma)_m}{m!}\frac{1}{(\frac{3}{\delta}+m)(\frac{2}{\delta}+m)}
\times \nonumber\\
& &
\left[\frac{\gamma!}{(\frac{2}{\delta}+m+1)\ldots(\frac{2}{\delta}
+m+\gamma)}\right. \nonumber \\
& & \left.-(\frac{r}{R})^{m\delta+2}\;_2F_1(-\gamma,
\frac{2}{\delta}+m;
\frac{2}{\delta}+m+1;(\frac{r}{R})^\delta)\right].\label{line3}
\end{eqnarray}
Then
\begin{equation}
P(r)=4\pi G\rho_c^2R^2g(r)
\end{equation}
and the temperature is then given by
\begin{equation}
T(r)=\frac{\mu}{kN_A}4\pi G\rho_cR^2\frac{g(r)}{[1-
(\frac{r}{R})^\delta]^\gamma}.
\end{equation}
In the following table we tabulate $g(r)$ for $\delta=1.28$ and
$\gamma=10$ from which we can compute the pressure for  the solar
core of the sum ranging from $0\leq \frac{r}{R}\leq 0.3.$\par
\bigskip
\noindent
\begin{tabbing}
\=Table 4 \=pressure = const.$\times g(r)$, temperature =
const.$\times g(r)/u,$\\
\>\>$u=\left[1-
(\frac{r}{R})^\delta\right]^\gamma$, for $\delta=1.28$ and
$\gamma=10$
\end{tabbing}
\begin{tabular}{lccc}
$y=\frac{r}{R}$ & $g(r)$ & $v=\left[1-
(\frac{r}{R})^\delta\right]^\gamma$ & $\frac{g(r)}{u}$\\ [0,3cm] 
0.0864 & $1.7911\times10^{-3}$ & 0.6418 & $2.7907\times10^{-3}$\\
0.1153 & $1.4035\times10^{-3}$ & 0.5228 & $2.6842\times10^{-3}$\\
0.1441 & $1.0546\times10^{-3}$ & 0.4179 & $2.5234\times10^{-3}$\\
0.1873 & $6.4117\times10^{-4}$ & 0.2884 & $2.2227\times10^{-3}$\\
0.2161 & $4.4168\times10^{-4}$ & 0.2202 & $2.0055\times10^{-3}$\\
0.2450 & $2.9482\times10^{-4}$ & 0.1649 & $1.7876\times10^{-3}$\\
0.2882 & $1.5279\times10^{-4}$ & 0.1033 & $1.4788\times10^{-3}$\\
\end{tabular}
\clearpage
\noindent
Hence the proportional decrease in temperature is: \par
\medskip
\noindent
Table 5 temperature distribution\\
\begin{tabbing}
\=0000000000000000 \=0000000000000000 \=00000000000000 \= \kill
\>Sears (1964) \>$\frac{g(r)}{\left[\rho(r)/\rho_c\right]^{1/2}}$
\>$\frac{g(r)}{\left[\rho(r)/\rho_c\right]^{1/4}}$\\ [0,3cm]
\>0.8789 \>0.9428 \>0.8987\\
\>0.9275 \>0.9316 \>0.8819\\
\>0.8828 \>0.9219 \>0.8385\\
\>0.8938 \>0.9612 \>0.8956\\
\>0.8911 \>0.9709 \>0.8992\\
\>0.9000 \>0.7644 \>0.8656
\end{tabbing}
\noindent
Table 6 pressure = const.$\times g(r)$, for $\delta=1.28$ and
$\gamma=10$
\begin{tabbing}
\=0,0000111122223333 \=0,0000111122223333 \= \kill
\>$y=\frac{r}{R}$ \>$g(r)$ \\
\>$0.0864$ \>$3.1681\times10^{-3}$\\
\>$0.1153$ \>$2.7142\times 10^{-3}$\\
\>$0.1441$ \>$2.2664\times10^{-3}$\\
\>$0.1873$ \>$1.7286\times10^{-3}$\\
\>$0.2161$ \>$1.4425\times10^{-3}$\\
\>$0.2450$ \>$1.2012\times10^{-3}$\\
\>$0.2882$ \>$0.9129\times10^{-3}$\\
\end{tabbing}
\section{ENERGY GENERATION}
From Mathai and Haubold (1988) we have the net release of
thermonuclear
energy per gram per second given by
\begin{eqnarray}
\epsilon[\rho(r), T(r)] & = & \epsilon_0 [\rho (r)]^n[T(r)]^m
\nonumber \\
& = & \epsilon_0[\rho(r)]^n[\frac{\mu}{kN_A}\frac{P(r)}{\rho(r)}]^m
\nonumber \\
& = & \epsilon_0(\frac{\mu}{kN_A})^m[\rho(r)]^{n-m}[P(r)]^m
\label{line3}
\end{eqnarray}
for the parameter n and m determined by the proton-proton chain and
the CNO cycle. Since the g(r) of table 5 is proportional to P(r) we
have
\begin{equation}
\epsilon[\rho(r), T(r)]=k_1[\rho(r)]^{n-m}[g(r)]^m,
\end{equation}
where $k_1$ is a constant.
The equation of energy conservation is
\begin{equation}
\frac{dL(r)}{dr}=4\pi r^2\rho(r)\epsilon[\rho(r), T(r)],
\end{equation}
where L(r) denotes the luminosity of the Sun. Hence the total
luminosity is given by
\begin{eqnarray}
L & = & 4\pi \int^R_0dr r^2 \rho(r)\epsilon[\rho(r), T(r)]\\
& = & 4\pi\epsilon_0\left(\frac{\mu}{kN_A}\right)^m \int^R_0dr
r^2[\rho(r)]^{1+n-m}[P(r)]^m
\end{eqnarray}
for the model in (3.2). Substituting for P(r) from (2.10)
one has
\begin{eqnarray}
L & = & (4\pi\epsilon_0)^{m+1}(\frac{GR^2}{\delta^2})^m
\rho_c^{m+n+1}(\frac{\mu}
{kN_A})^m \times \nonumber \\
& &
\int^R_0drr^2[1-(\frac{r}{R})^\delta]^{\gamma(1+n-m)}[\Psi(\gamma,
\delta)-\phi(\frac{r}{R})]^m, \label{line2}
\end{eqnarray}
where
\begin{eqnarray}
\phi(\frac{r}{R}) & = & \sum^\gamma_{m_1=0}\frac{(-
\gamma)_{m_1}}{m_1!}\frac{(\frac{r}{R})^{m_1\delta+2}}
{(\frac{3}{\delta}+m_1)(\frac{2}{\delta}+m_1)} \times \nonumber \\
& &_2F_1(-\gamma,
\frac{2}{\delta}+m_1;\frac{2}{\delta}+1+m_1;(\frac{r}{R}^\delta))
\label{line2}
\end{eqnarray}
and
\begin{equation}
\psi(\delta, \gamma)=\sum^\gamma_{m_1=0}\frac{(-
\gamma)_{m_1}}{m_1!}\frac{1}{(\frac{3}{\delta}+m_1)(\frac{2}
{\delta}
+m_1)}\frac{\gamma!}{(\frac{2}{\delta}+m_1+1)\cdots(\frac{2}
{\delta}+m_1+\gamma)}.
\end{equation}
Note that $|\frac{\phi(\frac{r}{R})}{\psi(\delta,\gamma)}|\leq
1.$
If 
m and n are not assumed to be positive integers then the general
binomial expansion can be used but one has to check the convergence
conditions in the final sum. Assuming that m and n are positive
integers we get polynomials from the two factors in the integral in
(3.6). For a general positive integer m we get higher powers of
$\Phi(u)$ but $\Phi(u)$ itself is a double sum. Hence the
expression becomes complicated. For m=1 and n a positive integer
one can get some simpler representations. Hence
we consider this case. Then we have
\begin{equation}
L=(4\pi\epsilon_0)^2(\frac{GR^2}{\delta^2})\rho_c^{n+2}(\frac{\mu
}{kN_A})R^
3[\Phi(\delta, \gamma)I_0-I_1],
\end{equation}
where
\begin{equation}
I_0=\int^1_0 du u^2[1-
u^\delta]^{n\gamma}=\frac{(n\gamma)!}{\delta}\frac{\Gamma(\frac{3}
{\delta})}{\Gamma(\frac{3}{\delta}+n\gamma+1)}
\end{equation}
and
\begin{eqnarray}
I_1 & = & \int^1_0du u^2[1-u^\delta]^{n\gamma} \phi(u)\nonumber \\
& = &
\frac{\delta^2}{6}\sum^\gamma_{m_1=0}\sum^\gamma_{m_2=0}\frac{(-
\gamma)_{m_1}}{m_1!}\frac{(-
\gamma)_{m_2}}{m_2!}\frac{(\frac{2}{\delta})_{m_1}(\frac{3}
{\delta})_{m_1}}{(\frac{2}{\delta}+1)_{m_1}(\frac{3}{\delta}+1)_{
m_1}}
\times
\nonumber \\
& = &
\frac{(\frac{2}{\delta})_{m_1+m_2}}{(\frac{2}{\delta}+1)_{m_1+m_2}}
\int^1_0du u^2[1-u^\delta]^{n\gamma} u^{m_1\delta+m_2\delta +2}.
\label{line3}
\end{eqnarray}
But
\begin{eqnarray}
\int^1_0du u^2[1-u^\delta]^{n\gamma} u^
{m_1\delta+m_2\delta+2} & = & \frac{(n\gamma)!}{\delta}\frac{\Gamma
(\frac{5}{\delta}+ m_1+
m_2)}{\Gamma(\frac{5}{\delta}+n\gamma+1+m_1+m_2)}\nonumber \\
& = &
\frac{(n\gamma)!}{\delta}\frac{\Gamma(\frac{5}{\delta})}{\Gamma
(\frac{5}{\delta}+n\gamma+1)}\times \nonumber \\
& & \frac{(\frac{5}{\delta})_{m_1+m_2}}
{(\frac{5}{\delta}+n\gamma+1)_{m_1+m_2}}.\label{line3}
\end{eqnarray}
Substituting (3.12) in (3.11) and writing the resulting expansions
as Kamp\'{e} de F\'{e}rier's function we have
\begin{eqnarray}
I_1=\frac{\delta}{6} & (n\gamma)! & \frac{\Gamma(\frac{5}{\delta})}
{\Gamma(\frac{5}{\delta}+n\gamma+1)}\times \nonumber \\
& & F^{2:3:1}_{2:2:0}\left[(^1_1)\left|^{\frac{2}{\delta},
\frac{5}{\delta}:-\gamma, \frac{2}{\delta}, \frac{3}{\delta}:
-\gamma}_{\frac{2}{\delta}+1, \frac{5}{\delta}+n\gamma+1:
\frac{2}{\delta}+1, \frac{9}{\delta}+1:}\right. \right].
\end{eqnarray}
The luminosity at any given point r is available from (3.3). That
is,
\begin{eqnarray}
L(r)= & (4\pi\epsilon_0)^{m+1} & (\frac{GR^2}{\delta^2})^m
\rho_c^{m+n+1}
(\frac{\mu}{kN_A})^m \times\label{line1} \\
& & \int^r_0dt t^2\left[1-(\frac{t}{R})^\delta\right]^{\gamma(1+n-
m)}\left[\Psi(\gamma, \delta)-\Phi(\frac{t}{R})\right]^m.\nonumber
\end{eqnarray}
We will evaluate this for m=1. Then the two integrals to be
evaluated are
$$I_2=\int^r_0dt t^2\left[1-(\frac{t}{R})^\delta\right]^{n\gamma} $$
and
$$I_3=\int^r_0dt t^2\left[1-
(\frac{t}{R})^\delta\right]^{n\gamma}\Phi(\frac{t}{R}).$$
That is,
\begin{eqnarray}
I_2 & = & R^3\int^{r/R}_0 du u^2 [1-u^\delta]^{n\gamma}\nonumber \\
& = & R^3\sum^{n\gamma}_{\alpha=0}\frac
{(-n\gamma)_\alpha}{\alpha!}\int^{r/R}_0 du
u^{2+\alpha\delta}\nonumber \\
& = & R^3\sum^{n\gamma}_{\alpha=0}\frac{(-
n\gamma)_\alpha}{\alpha!}\frac{(r/{R})^{3+\alpha\delta}}{3+
\alpha\delta}\nonumber \\
& = & \frac{r^3}{3}\sum^{n\gamma}_{\alpha=0}\frac{(-
n\gamma)_\alpha}{\alpha!}\frac{(\frac{3}{\delta})_\alpha}{(\frac{
3}{\delta}+1)_\alpha}\left[\left(\frac{r}{R}\right)^\delta\right]
^\alpha \nonumber \\
& = & \frac{r^3}{3}\;_2F_1\left(-n\gamma,
\frac{3}{\delta};\frac{3}{\delta}+1;\left(\frac{r}{R}\right)^
\delta\right).
\end{eqnarray}
\begin{eqnarray}
I_3 & = &
R^3\sum^\gamma_{m_1=0}\sum^\gamma_{m_2=0}\sum^{n\gamma}_{m_3=0}
\frac{(-\gamma)_{m_1}}{m_1!}\frac{(-
\gamma)_{m_2}}{m_2!}\frac{(-n\gamma)_{m_3}}{m_3!}\times
\label{line1}\\
& &
\frac{\left(\frac{2}{\delta}\right)_{m_1+m_2}}{\left(\frac{2}
{\delta}+1\right)_{m_1+m_2}}
\frac{\left(\frac{r}{R}\right)^{5+(m_1+m_2+m_3)\delta}}{\delta
\left[\frac{5}
{\delta}+m_1+m_2+m_3\right]}\nonumber \\
& = &
\frac{r^5}{5R^2}\sum^\gamma_{m_1=0}\sum^\gamma_{m_2=0}\frac{(-
\gamma)_{m_1}}{m_1!}\frac{(-
\gamma)_{m_2}}{m_2!}\frac{(\frac{5}{\delta})_{m_1+m_2}}{(\frac{5}
{\delta}+1)_{m_1+m_2}} \times \nonumber \\
& &
\frac{(\frac{2}{\delta})_{m_1+m_2}}{(\frac{2}{\delta}+1)_{m_1+m_2
}}\left[\left(\frac{r}{R}\right)^\delta\right]^{m_1+m_2}\;\times
\nonumber \\
& & _2F_1 \left(-n\gamma,
\frac{5}{\delta}+m_1+m_2;\frac{5}{\delta}+1+m_1+m_2;\left(\frac{r}
{R}\right)^\delta\right).\nonumber
\end{eqnarray}
Note that (3.16) can also be written in terms of a Kamp\'{e} de
F\'{e}riet's function of three variables. In that case the
variables
are $1, (r/R)^\delta, (r/R)^\delta,$ respectively. Note also that
(3.15) and (3.16) are polynomials. Hence for m=1, n a positive
integer
\begin{equation}
L(r)=(4\pi)^2\left(\frac{GR^2}{\delta^2}\right)\rho^{n+2}_c\left(
\frac{\mu}{kN_A}\right)\left[\Psi(\gamma, \delta)I_2-I_3\right],
\end{equation}
where $I_2$ and $I_3$ are given in (3.15) and (3.16) respectively.
A particular case of (3.17) is also available from Haubold and
Mathai (1992).
\vspace{2cm}
\begin{center}
{\bf REFERENCES}
\end{center}
H.J. Haubold and A.M. Mathai, Astrophys. Space Sci.
\underline{197}, 153(1982).\par
\medskip
\noindent
A.M. Mathai and H.J. Haubold, Modern Problems in Nuclear and
Neutrino Astrophysics. Akademie-Verlag, Berlin 1988.\par
\medskip
\noindent
A.M. Mathai, A Handbook of Generalized Special Functions for
Statistical and Physical Sciences. Oxford University Press, Oxford
1993.\par
\medskip
\noindent
A. Noels, R. Papy, and F. Remy, Computers in Physics \underline{7},
22(1993).\par
\medskip
\noindent
R.L. Sears, Astrophys. J. \underline{140}, 477 (1964).\par
\medskip
\noindent
H.M. Srivastava and P.W. Karlsson, Multiple Gaussian Hypergeometric
Series. Ellis Horwood, Chichester 1985.
\end{document}